\documentclass[a4paper,12pt]{article}
\usepackage{times, url}
\usepackage{hyperref}

\usepackage{graphicx,amsmath,amssymb,amsthm,amsfonts,bm,float,cite}
\newtheorem{theorem}{Theorem}[section]

\newtheorem{lemma}{Lemma}[section]

\usepackage[none]{hyphenat}

\usepackage[center]{caption}
\renewcommand{\thefootnote}{\fnsymbol{footnote}}

\usepackage{color, colortbl}
\definecolor{Gray}{gray}{0.9}

\newcommand\blfootnote[1]{%
  \begingroup
  \renewcommand\thefootnote{}\footnote{#1}%
  \addtocounter{footnote}{-1}%
  \endgroup
}

\usepackage[bottom,symbol]{footmisc}
\DeclareFontEncoding{TS1}{}{}
\DeclareFontSubstitution{TS1}{cmr}{m}{n}
\DeclareTextSymbol{\tcrp}{TS1}{'251}
\DeclareTextSymbolDefault{\tcrp}{TS1}
\newcommand{\licen}{ \phantom{.} \\  \hspace{-10mm} {\small  \fbox{
	\begin{tabular}{p{0.105\textwidth} p{0.82\textwidth}}
	  \raisebox{-22pt}{\includegraphics[scale=0.34]{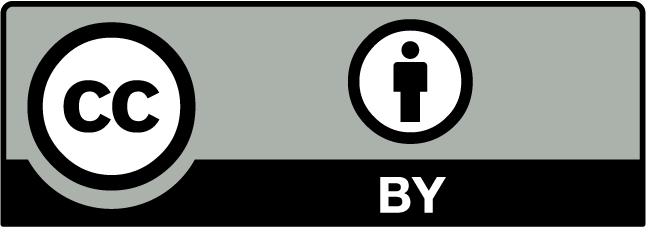}}	& {
			Copyright \tcrp ~2024 by the Authors. This is an Open Access paper distributed under the terms and conditions of the Creative Commons Attribution 4.0 International License (CC BY 4.0).} { {\footnotesize \ttfamily https://creativecommons.org/licenses/by/4.0/}}    \vspace{-0.5mm}
	\end{tabular} }
} \vspace{-10mm}
}

\begin{document}
\setcounter{page}{1}

\thispagestyle{empty} 

\noindent \textbf{Notes on Number Theory and Discrete
Mathematics \newline
Print ISSN 1310--5132, Online ISSN 2367--8275 \newline
XXXX, Volume XX, Number X, XXX--XXX \newline
DOI: 10.7546/nntdm.XXXX}  
\vspace{11mm}

\begin{center}
{\LARGE \bf  Elementary proof of W. Schulte's conjecture} 
\vspace{8mm}

{\Large \bf Djamel Himane$^1$}
\vspace{3mm}

$^1$ LA3C Laboratory, Faculty of Mathematics, USTHB, \\ 
Po. Box 32, El Alia, 16111, Bab Ezzouar, Algiers, Algeria. \\
e-mail: \url{dhimane@usthb.dz}
\vspace{2mm}
\end{center}
\vspace{9mm}  \blfootnote{\licen}

\noindent
\noindent {\bf Received:}  DD Month XXXX \hfill {\bf Revised:} DD Month XXXX \\
{\bf Accepted:} DD Month XXXX  \hfill {\bf Online First:} DD Month XXXX \\[4mm] 

{\bf Abstract:} In the On-Line Encyclopedia of Integer Sequences, we find the sequence \( 1, 1, 3, 18, \\ 180, 2700, 56700, 1587600, 57153600, \dots \) , which is given by the formula
$ A_{n} = n!(n-1)!/2^{n-1}. $

On the same page, Werner Schulte conjectured that for all \( n > 1 \), \( n \) divides \( 2A_{n-1} + 4 \) if and only if \( n \) is prime. In this paper, we employ elementary methods to provide a simple proof of this conjecture.
\vspace{5mm}

\section{Introduction and statement of results}
Consider the sequence defined by the relation:
\[
A_n = 1 \cdot (1 + 2) \cdot (1 + 2 + 3) \cdots (1 + 2 + \ldots + (n-1)).
\]
We pose the question: does this sequence of positive integers hold any particular significance?

In the course of analyzing this sequence, we consider the triangular number:
\[
T_s = 1 + 2 + \ldots + s = \frac{s(s + 1)}{2}.
\]
Using this formula, we can express the sequence \( A_n \) as follows:
\begin{equation} \label{eq.1} 
A_n = \prod_{j=1}^{n-1} T_j = \prod_{j=1}^{n-1} \frac{j(j + 1)}{2} = \frac{(n - 1)! \cdot n!}{2^{n-1}}.
\end{equation}

The first terms of this sequence are \( 1, 1, 3, 18, 180, 2700, 56700, \ldots \) (sequence \href{https://oeis.org/A006472}{A006472} in the OEIS, also listed as sequence M3052 in Simon \cite{1}).

In the OEIS database entry \href{https://oeis.org/A006472/internal}{A006472}, we encountered a conjecture attributed to W. Schulte regarding this sequence, for which we will provide an elementary proof.

The main result we aim to prove is the following:
\begin{theorem}[W. Schulte's Conjecture]
For any positive integer \( n > 1 \), \( n \) divides
\begin{equation} \label{eq.23} 
    \frac{(n-2)!(n-1)!}{2^{n-3}} + 4.
    \end{equation}
if and only if \( n \) is a prime number.
\end{theorem}

\section{Preliminaries}
For any prime \( p \) and any non-negative integer \( n \), we define \( v_{p}(n) \) as the exponent of the largest power of \( p \) that divides \( n \) (i.e., the \( p \)-adic valuation of \( n \)). We can express this using the following formula:
\begin{equation} \label{eq.2} 
    v_{p}(n!) = \sum_{i=1}^{\infty} \left\lfloor \frac{n}{p^{i}} \right\rfloor,
\end{equation}
where \( \lfloor x \rfloor \) denotes the floor function of \( x \).

Legendre’s formula \eqref{eq.2} provides an expression for the largest power of a prime \( p \) that divides the factorial \( n! \) \cite{2}. Although this formula is an infinite sum, for each pair of values \( n \) and \( p \), only a few terms are non-zero: for each \( i \) such that \( p^{i} > n \), we have \( \left\lfloor \dfrac{n}{p^{i}} \right\rfloor = 0 \).

We also have the following recurrence relation:
\begin{equation}
    v_{p}(n!) = \left\lfloor \dfrac{n}{p} \right\rfloor + v_{p}\left(\left\lfloor \dfrac{n}{p} \right\rfloor!\right).
\end{equation}
Legendre’s formula can also be rewritten in terms of the base-\( p \) expansion of \( n \). If \( s_{p}(n) \) denotes the sum of the digits in the base-\( p \) representation of \( n \), then we have:
\begin{equation}
    v_{p}(n!) = \dfrac{n - s_{p}(n)}{p-1}.
\end{equation}
It is easy to verify the following relations: For any positive integers \( a, b \) and any prime \( p \) \cite{3}:
\begin{enumerate}
    \item \( v_{p}(a \times b) = v_{p}(a) + v_{p}(b). \)
    \item \( v_{p}\left(\dfrac{a}{b}\right) = v_{p}(a) - v_{p}(b). \)
    \item If \( a \mid b \), then \( v_{p}(a) \leq v_{p}(b) \), with \( v_{p}(a) < a \).
\end{enumerate}

The proof of the following lemma is well-known and straightforward, but we include it here for completeness.

\begin{lemma} \label{Lemma 1} 
For any positive integer \( m > 4 \), we have \( m < 2^{m/2} \).
\end{lemma}

\begin{proof} 
We prove this inequality by induction. Note that in the induction step, the left side grows additively, whereas the right side grows multiplicatively. In particular, we have:
\begin{equation*}
    m < 2^{m/2},
\end{equation*} 
which is equivalent to
\begin{equation*}  
    m \left( \frac{m+1}{m} \right) < 2^{m/2} \left(\frac{m+1}{m}\right).
\end{equation*}

Simplifying, we obtain:
\begin{equation} \label{eq.3} 
    m + 1 < 2^{m/2} \left(\frac{m+1}{m}\right).
\end{equation}

We define the following auxiliary function:
\begin{equation*}
    \mathrm{g}(t) = \sqrt{2} - \frac{t+1}{t} = (\sqrt{2} - 1) - \frac{1}{t}.
\end{equation*}

Evaluating, we find \( \mathrm{g}(5) = (\sqrt{2} - 1) - \frac{1}{5} \approx 0.2142 \). Moreover, since \( \mathrm{g}'(t) = \frac{1}{t^{2}} > 0 \), we see that \( \mathrm{g}(t) \) is a strictly increasing function for \( t \geq 5 \). Thus, for \( t \geq 5 \), we have
\begin{equation*}
    \sqrt{2} - \frac{t+1}{t} > 0 \quad \Rightarrow \quad \sqrt{2} > \frac{t+1}{t}.
\end{equation*}

This implies that, for \( m \geq 5 \),
\begin{equation*}
    \frac{m+1}{m} < 2^{1/2} \quad \Rightarrow \quad 2^{m/2} \left( \frac{m+1}{m} \right) < 2^{(m+1)/2}.
\end{equation*}

Using this relation in \eqref{eq.3}, we obtain:
\begin{equation*}
    m + 1 < 2^{(m+1)/2},
\end{equation*}
which completes the induction step. Therefore, for \( m > 4 \), it follows that
\begin{equation*}
    m < 2^{m/2}.
\end{equation*}
\end{proof}

\begin{lemma} \label{Lemma 2}
If \( p \) is prime, then \( (p-1)! + 1 \) is divisible by \( p \).
\end{lemma}

This theorem was proposed by Alhazen (c. 965-1000 AD) and, in the 18th century, by J. Wilson. Lagrange provided the first proof in 1771. A primality test \cite{4} directly follows from Wilson’s theorem: a positive integer \( n > 1 \) is prime if and only if 
\begin{equation}
    (n-1)! + 1 \equiv 0 \pmod{n}.
\end{equation}
  
\begin{lemma}  \label{Lemma 3}
For any positive integer \( a \) that is not divisible by a prime \( p \), we have
\begin{equation}
    a^{p-1} \equiv 1 \pmod{p}.
\end{equation}
\end{lemma}

This fundamental result in number theory is known as Fermat's Little Theorem. First stated by Pierre de Fermat in 1640, this theorem forms the basis for Fermat’s primality test  (p. 22, \cite{4}).

\section{Proof of the main theorem}

\begin{proof}
Let \( n \) be a positive integer greater than 1. We will prove that \( n \) divides
\[
\dfrac{(n-2)! \, (n-1)!}{2^{n-3}} + 4
\]
if and only if \( n \) is a prime number. The proof is divided into two parts:

\begin{itemize}
    \item[$\bullet$] \textbf{First part:} Assume that \( n \) is a prime number. Let \( n = p \), where \( p \) is an odd prime. By Wilson's theorem (\ref{Lemma 2}), we know that \( p \) divides both \( (p-1)! + 1 \) and \( (p-2)! - 1 \). Moreover, by Fermat's Little Theorem (\ref{Lemma 3}), we deduce that \( p \) divides \( 2^{p-1} - 1 \). Using the relation \eqref{eq.1}, we obtain:
    \[ 2 \times A_{p-1} \equiv -4 \pmod{p}. \]
    Therefore, \( p \) divides
    \[
    \dfrac{(p-2)! \, (p-1)!}{2^{p-3}} + 4.
    \]
    For \( n = 2 \), the property holds directly.

    \item[$\bullet$] \textbf{Second part:} Now assume that \( n \) is a composite number. Let \( n = pq \), where \( 1 < p \leq q < n \). We consider three cases:
    \begin{itemize}
        \item[$\star$] \textbf{Case 1:} If \( n = 4 \), we have
        \[
        \dfrac{(n-2)! \, (n-1)!}{2^{n-3}} + 4 = 10,
        \]
        which is not divisible by 4.

        \item[$\star$] \textbf{Case 2:} If \( n \) is odd and \( n > 4 \), then \( n \mid (n-2)! \, (n-1)! \) and \( \gcd(n, 2^{n-3}) = 1 \). We deduce that:
        \[
        n \mid \dfrac{(n-2)! \, (n-1)!}{2^{n-3}},
        \]
        which implies:
        \[
        \dfrac{(n-2)! \, (n-1)!}{2^{n-3}} + 4 \equiv 4 \pmod{n}.
        \]
        This means that the expression is not divisible by \( n \).

        \item[$\star$] \textbf{Case 3:} If \( n \) is even, let \( n = 2^{v_{2}(n)} \times k > 4 \), where \( k \) is odd. We have \( n \mid (n-2)! (n-1)! \). By Euclid's lemma, we know that \( k \mid 2 \times A_{n-1} \), which implies that \( n \mid 2 \times A_{n-1} \) if and only if
        \[
        v_{2}(n) \leq v_{2}(2 \times A_{n-1}).
        \]
        Using properties of the \( p \)-adic valuation, we obtain:
        \[
        v_{2}(2 \times A_{n-1}) = v_{2}((n-1)!) + v_{2}((n-2)!) - v_{2}(2^{n-3}).
        \]
        We use the recurrence relation:
        \[
        v_{2}((n-1)!) = v_{2}(n-1) + v_{2}((n-2)!) = v_{2}((n-2)!),
        \]
        since \( n-1 \) is odd, and thus
        \[
        v_{2}((n-2)!) = \left\lfloor \dfrac{n-2}{2} \right\rfloor + \left\lfloor \dfrac{n-2}{4} \right\rfloor + v_{2}\left(\left\lfloor \dfrac{n-2}{4} \right\rfloor!\right) \geq \dfrac{n-2}{2} + \dfrac{n-2}{4}.
        \]
        Therefore,
        \[
        v_{2}(2 \times A_{n-1}) \geq \left( \dfrac{n-2}{2} + \dfrac{n-2}{4} \right) + \left( \dfrac{n-2}{2} + \dfrac{n-2}{4} \right) - (n-3).
        \]
        Simplifying, we find:
        \[
        v_{2}(2 \times A_{n-1}) \geq \dfrac{n}{2}.
        \]
        By Lemma \ref{Lemma 1}, this implies that
        \[
        \dfrac{n}{2} > v_{2}(n).
        \]
        Thus,
        \[
        v_{2}(2 \times A_{n-1}) \geq v_{2}(n),
        \]
        which implies that
        \[
        \dfrac{(n-2)! \, (n-1)!}{2^{n-3}} + 4 \equiv 4 \pmod{n}.
        \]
        In conclusion,
        \[
        \dfrac{(n-2)! \, (n-1)!}{2^{n-3}} + 4
        \]
        is not divisible by \( n \).
    \end{itemize}
\end{itemize}

Thus, the proof is complete.
\end{proof}

\section*{Acknowledgements} 

 I would like to express my sincere appreciation to Professor Boumahdi Rachid at NHSM, who generously offered his advice and guidance. I would also like to extend my deepest gratitude to the reviewers for their voluntary efforts in supporting research and the scientific community.
 
\makeatletter
\renewcommand{\@biblabel}[1]{[#1]\hfill}
\makeatother

\end{document}